\documentclass[11pt]{article}

\usepackage{amssymb}
\usepackage{amsfonts}
\usepackage{amsmath}
\usepackage{amsthm}
\usepackage[margin=1.00in]{geometry}
\usepackage{enumerate}
\usepackage{amstext}
\usepackage{layout}

\numberwithin{equation}{section}

\newtheorem{theorem}{Theorem}[section]
\newtheorem{corollary}{Corollary}[section]
\newtheorem{lemma}{Lemma}[section]

\newtheorem{remark}{Remark}[section]

\begin{document}

\begin{center}
{\large \bf A probability inequality for sums of independent Banach space valued random variables}

\vskip 0.5cm

{\bf Deli Li$^{1}$, 
Han-Ying Liang$^{2,}$\footnote{Corresponding author: Han-Ying Liang (Telephone: 86-21-65983242, FAX: 86-21-65983242)}, 
and Andrew Rosalsky$^{3}$}

\vskip 0.3cm

$^1$Department of Mathematical Sciences, Lakehead
University,\\
Thunder Bay, Ontario, Canada\\
$^2$Department of Mathematics, Tongji University,\\
Shanghai, China\\
$^3$Department of Statistics, University of Florida,\\
Gainesville, Florida, USA

\end{center}

\vskip 0.3cm

\noindent {\bf Abstract}~~Let $(\mathbf{B}, \|\cdot\|)$ be a real separable Banach space. 
Let $\varphi(\cdot)$ and $\psi(\cdot)$ be two continuous and increasing functions
defined on $[0, \infty)$ such that $\varphi(0) = \psi(0) = 0$, $\lim_{t \rightarrow \infty} \varphi(t) = \infty$, 
and $\frac{\psi(\cdot)}{\varphi(\cdot)}$ is a nondecreasing function on $[0, \infty)$.
Let $\{V_{n};~n \geq 1 \}$ be a sequence of independent and symmetric {\bf B}-valued random variables.
In this note, we establish a probability inequality for sums of independent {\bf B}-valued 
random variables by showing that for every $n \geq 1$ and all $t \geq 0$,
\[
\mathbb{P}\left(\left\|\sum_{i=1}^{n} V_{i} \right\| > t b_{n} \right) 
\leq 4 \mathbb{P} \left(\left\|\sum_{i=1}^{n} \varphi\left(\psi^{-1}(\|V_{i}\|)\right) 
\frac{V_{i}}{\|V_{i}\|} \right\| > t a_{n} \right) 
+ \sum_{i=1}^{n}\mathbb{P}\left(\|V_{i}\| > b_{n} \right),
\]
where $a_{n} = \varphi(n)$ and $b_{n} = \psi(n)$, $n \geq 1$. As an application of this inequality, we establish 
what we call a comparison theorem for the weak law of large numbers for independent and identically distributed 
${\bf B}$-valued random variables.

~\\


\vskip 0.3cm


\vskip 0.3cm

\noindent {\bf Running Head}: Probability inequality

\section{Introduction and the main result}

We begin with some notation. Let $(\Omega, \mathcal{F}, \mathbb{P})$ be a probability space and
let $(\mathbf{B}, \| \cdot \| )$ be a real separable Banach space 
equipped with its Borel $\sigma$-algebra $\mathcal{B}$ 
($=$ the $\sigma$-algebra generated by the class of open subsets of
$\mathbf{B}$ determined by $\|\cdot\|$). A {\bf B}-valued random variable
$X$ is defined as a measurable function from $(\Omega, \mathcal{F})$ into $(\mathbf {B}, \mathcal{B})$.
Throughout this note, let $\{R, R_{n}; ~n \geq 1\}$ be a Rademacher sequence; i.e., $\{R, R_{n}; ~n \geq 1\}$ 
is a sequence of independent and identically distributed (i.i.d.) real-valued random variables with 
$\mathbb{P}(R = -1) = \mathbb{P}(R = 1) = 1/2$. 

Probability inequalities are essential for establishing virtually all
probability limit theorems and for advancing statistical theory, and
they also have an intrinsic interest as well. For example, some invaluable and
celebrated classical inequalities are those of L\'{e}vy, Ottaviani, Kahane, 
Hoffmann-J{\o}rgensen (see, e.g., L\'{e}vy (1937), Chow and
Teicher (1997, p. 75), Kahane (1968), and Hoffmann-J{\o}rgensen (1974),
respectively), etc. This note is devoted to establishing the following 
probability inequality which is a comparison theorem for sums of independent 
$\mathbf{B}$--valued random variables.

\vskip 0.3cm

\begin{theorem}
Let $\varphi(\cdot)$ and $\psi(\cdot)$ be two continuous and increasing functions
defined on $[0, \infty)$ such that $\varphi(0) = \psi(0) = 0$ and
\begin{equation}
\lim_{t \rightarrow \infty} \varphi(t) = \infty 
~~\mbox{and}~~\frac{\psi(\cdot)}{\varphi(\cdot)} ~\mbox{is a nondecreasing function on}~[0, \infty).
\end{equation}
Here we define $\frac{\varphi(0)}{\psi(0)} = \lim_{t \rightarrow 0^{+}} \frac{\varphi(t)}{\psi(t)}$.
For $n \geq 1$, set $a_{n} = \varphi(n)$ and $b_{n} = \psi(n)$. Then we have:

{\bf (i)}~~Let $\{x_{n};~ n \geq 1\}$ be a {\bf B}-valued sequence such that $\|x_{n}\| \leq b_{n}, ~n \geq 1$. 
Then for every $n \geq 1$ and all $t \geq 0$,
\begin{equation}
\mathbb{P}\left(\left\|\sum_{i=1}^{n} R_{i}x_{i} \right\| > t b_{n} \right) 
\leq 2 \mathbb{P} \left(\left\|\sum_{i=1}^{n} R_{i}
\varphi\left(\psi^{-1}(\|x_{i}\|)\right) \frac{x_{i}}{\|x_{i}\|} \right\| > t a_{n} \right).
\end{equation}
Here $\displaystyle \varphi\left(\psi^{-1}(\|0\|)\right) \frac{0}{\|0\|} \stackrel{\Delta}{=} 0$ since 
$\displaystyle \lim_{x \rightarrow 0} \varphi\left(\psi^{-1}(\|x\|)\right) \frac{x}{\|x\|} = 0$.

{\bf (ii)}~~If $\{V_{n};~n \geq 1 \}$ is a sequence of independent 
and symmetric {\bf B}-valued random variables, then for every 
$n \geq 1$ and all $t \geq 0$,
\begin{equation}
\mathbb{P}\left(\left\|\sum_{i=1}^{n} V_{i} \right\| > t b_{n} \right) 
\leq 4 \mathbb{P} \left(\left\|\sum_{i=1}^{n} \varphi\left(\psi^{-1}(\|V_{i}\|)\right) 
\frac{V_{i}}{\|V_{i}\|} \right\| > t a_{n} \right) 
+ \sum_{i=1}^{n}\mathbb{P}\left(\|V_{i}\| > b_{n} \right).
\end{equation}
\end{theorem}

\vskip 0.3cm

The proof of Theorem 1.1 is given in Section 2. An application of Theorem 1.1 is presented in Section 3.
The application provides what we refer to as a comparison theorem for the weak law of large numbers (WLLN) 
for i.i.d. ${\bf B}$-valued random variables.

\section{Proof of Theorem 1.1}

To establish Theorem 1.1, we use the following result which is the second part of Theorem 4.4 of Ledoux and Talagrand (1991).  

\begin{lemma}
Let $\{x_{n};~ n \geq 1\}$ be a {\bf B}-valued sequence and let $\{\alpha_{n}; ~n \geq 1 \}$ be a real-valued
sequence such that $\sup_{n \geq 1}|\alpha_{n}| \leq 1$. Then we have, for every $n \geq 1$ and all $t \geq 0$,
\[
\mathbb{P} \left(\left\|\sum_{i=1}^{n} \alpha_{i} R_{i} x_{i} \right\| > t \right) 
\leq 2 \mathbb{P} \left(\left\|\sum_{i=1}^{n} R_{i} x_{i} \right\| > t \right).
\]
\end{lemma}

\vskip 0.3cm

\noindent {\it Proof of Theorem 1.1}~~Part {\bf (i)}. Clearly we have, for every 
$n \geq 1$ and all $t \geq 0$,
\[
\begin{array}{lll}
\mbox{$\displaystyle 
\mathbb{P}\left(\left\|\sum_{i=1}^{n} R_{i}x_{i} \right\| > t b_{n} \right)$}
&=& \mbox{$\displaystyle 
\mathbb{P}\left(\left\|\sum_{i=1}^{n} \left(\frac{a_{n}}{b_{n}} 
\cdot \frac{\psi\left(\psi^{-1}(\|x_{i}\|) \right)}{\varphi\left(\psi^{-1}(\|x_{i}\|) \right)}
\right) R_{i} \frac{\varphi\left(\psi^{-1}(\|x_{i}\|)\right)}{\psi\left(\psi^{-1}(\|x_{i}\|) \right)} x_{i}
\right\| > t a_{n} \right)$}\\
&&\\
&=& \mbox{$\displaystyle \mathbb{P}\left(\left\|\sum_{i=1}^{n} \left(\frac{\varphi(n)}{\psi(n)} 
\cdot \frac{\psi\left(\psi^{-1}(\|x_{i}\|) \right)}{\varphi\left(\psi^{-1}(\|x_{i}\|) \right)}
\right) R_{i} \varphi\left(\psi^{-1}(\|x_{i}\|)\right) \frac{x_{i}}{\|x_{i}\|}
\right\| > t a_{n} \right)$}.
\end{array}
\]
Since $\varphi(\cdot)$ and $\psi(\cdot)$ are two continuous and increasing functions defined on
$[0, \infty)$ satisfying $\varphi(0) = \psi(0) = 0$ and (1.1), we see that $\varphi^{-1}(\cdot)$ 
is also a continuous and increasing function defined on
$[0, \infty)$ such that $\psi^{-1}(0) = 0$, $\lim_{t \rightarrow \infty} \psi^{-1}(t) = \infty$, and
\[
0 \leq \frac{\psi\left(\psi^{-1}(t)\right)}{\varphi\left(\psi^{-1}(t)\right)} 
\leq \frac{\psi(n)}{\varphi(n)} = \frac{b_{n}}{a_{n}}~~\mbox{whenever}~~
0 \leq t \leq \psi(n) = b_{n}.
\]
Note that $\|x_{n}\| \leq b_{n} = \psi(n)$, $n \geq 1$.
We thus conclude that, for every $n \geq 1$,
\[
0 \leq \frac{\psi\left(\psi^{-1}(\|x_{i}\|) \right)}{\varphi\left(\psi^{-1}(\|x_{i}\|) \right)}
\leq \frac{\psi(n)}{\varphi(n)} ~~\mbox{for}~ i = 1, 2, ..., n
\]
and hence that, for every $n \geq 1$, 
\[
0 \leq \frac{\varphi(n)}{\psi(n)} 
\cdot \frac{\psi\left(\psi^{-1}(\|x_{i}\|) \right)}{\varphi\left(\psi^{-1}(\|x_{i}\|) \right)}
\leq 1~~\mbox{for}~ i = 1, 2, ..., n.
\]
By applying Lemma 2.1 we thus have, for every $n \geq 1$ and all $t \geq 0$,
\[
\begin{array}{ll}
& \mbox{$\displaystyle 
\mathbb{P}\left(\left\|\sum_{i=1}^{n} \left(\frac{\varphi(n)}{\psi(n)} 
\cdot \frac{\psi\left(\psi^{-1}(\|x_{i}\|) \right)}{\varphi\left(\psi^{-1}(\|x_{i}\|) \right)}
\right) R_{i} \varphi\left(\psi^{-1}(\|x_{i}\|)\right) \frac{x_{i}}{\|x_{i}\|}
\right\| > t a_{n} \right)$}\\
&\\
& \mbox{$\displaystyle 
\leq 2 \mathbb{P} \left(\left\|\sum_{i=1}^{n} R_{i}
\varphi\left(\psi^{-1}(\|x_{i}\|)\right) \frac{x_{i}}{\|x_{i}\|} \right\| > t a_{n} \right)$}
\end{array}
\]
proving Part {\bf (i)}.

Part {\bf (ii)}. For every $n \geq 1$, write
\[
V_{n,i} = V_{i}I\left\{\|V_{i} \| \leq b_{n} \right\},
~~T_{i} = \varphi\left(\psi^{-1}(\|V_{i}\|)\right)\frac{V_{i}}{\|V_{i}\|}, 
~~T_{n,i} = \varphi\left(\psi^{-1}(\|V_{n,i}\|)\right)\frac{V_{n,i}}{\|V_{n,i}\|}, ~~i = 1, ..., n.
\]
Clearly we have, for all $t \geq 0$,
\begin{equation}
\mathbb{P}\left(\left\|\sum_{i=1}^{n} V_{i} \right\| > t b_{n} \right) 
\leq \mathbb{P} \left( \left\|\sum_{i=1}^{n} V_{n,i} \right\| > t b_{n} \right) 
+ \sum_{i=1}^{n}\mathbb{P}\left(\|V_{i}\| > b_{n} \right).
\end{equation}
Note that
\[
\left\{\|V_{i} \| \leq b_{n} \right\} = 
\left\{\psi^{-1}(\|V_{i} \|) \leq n \right\} =
\left\{\varphi\left(\psi^{-1}(\|V_{i} \|) \right) \leq a_{n} \right\} 
= \left\{\|T_{i} \| \leq a_{n} \right\}, ~~i = 1, ..., n.
\]
Thus it is easy to see that
\[
T_{n,i} = T_{i} I\left\{\|T_{i} \| \leq a_{n} \right\}, ~~i = 1, ..., n.
\]
Since $\{V_{i};~i \geq 1 \}$ is a sequence of independent 
and symmetric {\bf B}-valued random variables,  $\{V_{n,i};~i = 1, ..., n \}$, $\{T_{i};~i = 1, ..., n \}$, and
$\{T_{n,i};~i = 1, ..., n \}$ are finite sequences of independent and symmetric {\bf B}-valued random variables.
Let $\{R, R_{n}; ~n \geq 1\}$ be a Rademacher sequence which is independent of $\{V_{n};~n \geq 1 \}$.
Then $\left\{R_{i}V_{n,i}; ~i = 1, ..., n \right\}$ has the same distribution as 
$\left\{V_{n,i}; ~i = 1, ..., n \right\}$ in ${\bf B}^{n}$ and 
$\left\{R_{i}T_{n,i};~i = 1, ..., n \right\}$ has the same distribution as 
$\left\{T_{n, i};~ i = 1, ..., n \right\}$ in ${\bf B}^{n}$. Since $\|V_{n,i}\| \leq b_{n}, ~i = 1, ..., n$, 
by applying (1.2), we have, for all $t \geq 0$,
\begin{equation}
\begin{array}{lll}
\mbox{$\displaystyle
\mathbb{P} \left( \left\|\sum_{i=1}^{n} V_{n,i} \right\| > t b_{n} \right)$}
&=& \mbox{$\displaystyle
\mathbb{P} \left( \left\|\sum_{i=1}^{n} R_{i} V_{n,i} \right\| > t b_{n} \right)$}\\
&&\\
&=& 
\mbox{$\displaystyle
\mathbb{E}\left(\mathbb{P}\left(\left. \left\|\sum_{i=1}^{n} R_{i} V_{n,i} \right\| > t b_{n} 
\right| V_{1}, ..., V_{n} \right) \right)$}\\
&&\\
&\leq&
\mbox{$\displaystyle 
2 \mathbb{E}\left(\mathbb{P}\left(\left. \left\|\sum_{i=1}^{n} R_{i} T_{n,i} \right\| > t a_{n} 
\right| V_{1}, ..., V_{n} \right) \right)$}\\
&&\\
&=&
\mbox{$\displaystyle 
2 \mathbb{P} \left( \left\|\sum_{i=1}^{n} R_{i} T_{n,i} \right\| > t a_{n} \right)$}\\
&&\\
&=&
\mbox{$\displaystyle 
2 \mathbb{P} \left( \left\|\sum_{i=1}^{n} T_{n,i} \right\| > t a_{n} \right)$}.
\end{array}
\end{equation}
For all $n \geq 1$, since $\{T_{i};~ i = 1, ..., n \}$ is a finite sequence of independent 
and symmetric {\bf B}-valued random variables, it follows that 
$\left\{T_{i}I\left\{\|T_{i}\| \leq a_{n} \right\} - T_{i}I\left\{\|T_{i}\| > a_{n} \right\}; ~i = 1, ..., n \right\}$ 
has the same distribution as $\{T_{i};~i = 1, ..., n \}$ in ${\bf B}^{n}$.
Note that 
\[
\sum_{i=1}^{n} T_{n,i} = \frac{\sum_{i=1}^{n}T_{i} + \sum_{i=1}^{n} \left(T_{i}I\left\{\|T_{i}\| \leq a_{n} \right\}
- T_{i}I\left\{\|T_{i}\| > a_{n} \right\} \right)}{2}, ~ n \geq 1.
\]
We thus have, for every $n \geq 1$ and all $t \geq 0$,
\begin{equation}
\begin{array}{lll}
\mbox{$\displaystyle
\mathbb{P} \left( \left\|\sum_{i=1}^{n} T_{n,i} \right\| > t a_{n} \right)$}
&\leq& \mbox{$\displaystyle
\mathbb{P} \left( \left\|\sum_{i=1}^{n} T_{i} \right\| > t a_{n} \right)$}\\
&&\\
&& \mbox{$\displaystyle
+ \mathbb{P} \left( \left\|\sum_{i=1}^{n} \left(T_{i}I\left\{\|T_{i}\| \leq a_{n} \right\}
- T_{i}I\left\{\|T_{i}\| > a_{n} \right\} \right) \right\| > t a_{n} \right)$}\\
&&\\
&=& \mbox{$\displaystyle
2 \mathbb{P} \left( \left\|\sum_{i=1}^{n} T_{i} \right\| > t a_{n} \right).$}
\end{array}
\end{equation}
Now we can see that (1.3) follows from (2.1), (2.2), and (2.3). ~$\Box$

\section{An application}

As an application of Theorem 1.1, in this section we will establish what we call 
a comparison theorem for the WLLN for i.i.d. {\bf B}-valued random variables. 
Theorem 3.1 is new even when $\mathbf{B} = \mathbb{R}$. 

\vskip 0.3cm

\begin{theorem}
{\rm{(Comparison theorem for the WLLN).}}
Let $(\mathbf{B}, \|\cdot\|)$ be a real separable Banach space.
Let $\{a_{n}; n \geq 1\}$ and $\{b_{n}; n \geq 1\}$ be increasing
sequences of positive real numbers such that 
\begin{equation}
\lim_{n \rightarrow \infty} a_{n} = \infty~ \mbox{and}~ \left\{b_{n}/a_{n}; ~n \geq 1 \right\}~
\mbox{is a nondecreasing sequence}.
\end{equation}
Suppose that, for every symmetric sequence $\{X, X_{n}; ~n \geq 1 \}$ of i.i.d. {\bf B}-valued random variables,
\begin{equation}
\frac{S_{n}}{a_{n}} \rightarrow_{\mathbb{P}} 0 ~~\mbox{if and only if}~~\lim_{n \rightarrow \infty}
n \mathbb{P}(\|X\| > a_{n}) = 0.
\end{equation}
Here and below $S_{n} = \sum_{i=1}^{n} X_{i},~n \geq 1$. Then, for every sequence 
$\{X, X_{n}; ~n \geq 1 \}$ of i.i.d. {\bf B}-valued random variables, we have that
\begin{equation}
\frac{S_{n}- \gamma_{n}}{b_{n}} 
\rightarrow_{\mathbb{P}} 0 ~~\mbox{or}~~\limsup_{n \rightarrow \infty} 
\mathbb{P} \left(\frac{\left\|S_{n}- \gamma_{n} \right\|}{b_{n}} 
> \lambda \right) > 0 ~~\forall~\lambda > 0
\end{equation}
according as
\begin{equation}
\lim_{n \rightarrow \infty}
n \mathbb{P}(\|X\| > b_{n}) = 0~~\mbox{or}~~\limsup_{n \rightarrow \infty}
n \mathbb{P}(\|X\| > b_{n}) >0.
\end{equation}
Here and below $\gamma_{n} = n \mathbb{E}\left(XI\{\|X\| \leq b_{n} \} \right)$, $n \geq 1$.
\end{theorem}

\vskip 0.2cm

\begin{remark}
Under the assumptions of Theorem 3.1, we conclude from Theorem 3.1 that, for every 
sequence $\{X, X_{n}; ~n \geq 1 \}$ of i.i.d. {\bf B}-valued random variables,
\[
\frac{S_{n} - \gamma_{n}}{b_{n}} \rightarrow_{\mathbb{P}} 0 ~~\mbox{if and only if}~~
\lim_{n \rightarrow \infty}
n \mathbb{P}(\|X\| > b_{n}) = 0,
\]
\[
\limsup_{n \rightarrow \infty} 
\mathbb{P} \left(\frac{\left\|S_{n}- \gamma_{n} \right\|}{b_{n}} 
> \lambda \right) > 0 ~~\forall~\lambda > 0 ~~\mbox{if and only if}~~
\limsup_{n \rightarrow \infty}
n \mathbb{P}(\|X\| > b_{n}) > 0.
\]
Hence 
\[
\frac{S_{n} - \gamma_{n}}{b_{n}} \nrightarrow_{\mathbb{P}} x ~~\forall ~x \in \mathbf{B}\backslash \{0\}.
\]
\end{remark}

\vskip 0.2cm

Let $0 < p \leq 2$. Then $\mathbf{B}$ is said to be of {\it stable type $p$} if
\[
\sum_{n=1}^{\infty} \Theta_{n}v_{n} ~~\mbox{converges a.s. whenever}~~
\{v_{n}: ~n \geq 1\} \subseteq \mathbf{B} ~~\mbox{with}~~
\sum_{n=1}^{\infty} \|v_{n}\|^{p} < \infty,
\]
where $\{\Theta_{n}; ~n \geq 1 \}$ is a sequence of i.i.d. stable random variables 
each with characteristic function $\psi(t) = \exp \left\{-|t|^{p}\right\}, ~- \infty < t <
\infty$. A remarkable characterization of stable type $p$ Banach spaces via the WLLN
was provided by Marcus and Woyczy\'{n}ski (1979) who showed that, for given $1 \leq p < 2$, the following two
statements are equivalent:
\begin{align*}
& {\bf (i)} \quad \mbox{The Banach space $\mathbf{B}$ is of
stable type $p$.}\\
& {\bf (ii)} \quad \mbox{For every symmetric sequence $\{X, X_{n}; ~n \geq 1 \}$
of i.i.d. {\bf B}-valued variables},
\end{align*}
\[
\frac{S_{n}}{n^{1/p}} \rightarrow_{\mathbb{P}} 0~~\mbox{if
and only if}~~\lim_{n \rightarrow \infty} n \mathbb{P}\left(\|X\| > n^{1/p}\right) = 0.
\]
Combining Theorem 3.1 and the above characterization of stable type $p$ Banach spaces, 
we immediately obtain the following two results.

\vskip 0.2cm

\begin{corollary}
Let $1 \leq p < 2$ and let $\{a_{n}; n \geq 1\}$ be an increasing sequence of positive 
real numbers such that 
\[
\lim_{n \rightarrow \infty} a_{n} = \infty~ \mbox{and}~ \left\{n^{1/p}/a_{n}; ~n \geq 1 \right\}~
\mbox{is a nondecreasing sequence}.
\]
Let $(\mathbf{B}, \|\cdot\|)$ be a real separable Banach space such that, for every symmetric sequence 
$\{X, X_{n}; ~n \geq 1 \}$ of i.i.d. {\bf B}-valued random variables, (3.2) holds. 
Then the Banach space {\bf B} is of stable type $p$. 
\end{corollary}

\vskip 0.2cm

\begin{corollary}
Let $(\mathbf{B}, \|\cdot\|)$ be a real separable Banach space.
Let $1 \leq p < 2$ and let $\{b_{n}; n \geq 1\}$ be a sequence of positive 
real numbers such that 
\[
\left\{b_{n}/n^{1/p}; ~n \geq 1 \right\}~\mbox{is a nondecreasing sequence}.
\]
If $\mathbf{B}$ is of stable type $p$, then for every sequence $\{X, X_{n}; ~n \geq 1 \}$ 
of i.i.d. {\bf B}-valued random variables, (3.3) and (3.4) are equivalent.
\end{corollary}

\vskip 0.2cm

\begin{corollary}
Let $\left \{X, X_{n};~ n \geq 1 \right \}$ be a sequence of i.i.d. real-valued random variables
and let $\left \{b_{n};~n \geq 1 \right \}$ be a sequence of positive real numbers such that $b_{n}/n^{1/p}$
is nondecreasing for some $p \in [1, 2)$. Set $S_{n} = \sum_{i=1}^{n} X_{i}, ~n \geq 1$. Then
\[
\frac{S_{n}- n \mathbb{E}(X I\{|X| \leq b_{n} \})}{b_{n}} 
\rightarrow_{\mathbb{P}} 0 ~~\mbox{or}~~\limsup_{n \rightarrow \infty} 
\mathbb{P} \left(\frac{\left|S_{n}- n \mathbb{E}(X I\{|X| \leq b_{n} \}) \right|}{b_{n}} 
> \lambda \right) > 0 ~~\forall~\lambda > 0
\]
according as
\[
\lim_{n \rightarrow \infty}
n \mathbb{P}(|X| > b_{n}) = 0~~\mbox{or}~~\limsup_{n \rightarrow \infty}
n \mathbb{P}(|X| > b_{n}) >0.
\]
\end{corollary}

\vskip 0.2cm

\noindent {\it Proof}.~~It is well known that the real line $\mathbb{R}$ is of stable type $p$ for all 
$p \in [1, 2)$ and so the corollary follows immediately from Corollary 3.2. ~$\Box$

\vskip 0.2cm

\begin{remark}
Corollary 3.3 is an improved version of Theorem 1 (ii) of Klass and Teicher (1977) wherein 
$\left \{b_{n}/n;~ n \geq 1 \right \}$ is nondecreasing. Theorem 1 of Klass and Teicher (1977)
may be regarded as a WLLN analogue of Feller's (1946) extension of Marcinkiewicz-Zygmund SLLN 
(see, e.g., Chow and Teicher (1997, p. 125)) to more general norming sequences. 
\end{remark}

\vskip 0.3cm

To prove Theorem 3.1, we use the following two preliminary lemmas. The second lemma is due to
Li, Liang, and Rosalsky (2016).

\vskip 0.2cm

\begin{lemma}
Let $\left \{a_{n};~n \geq 1 \right \}$ and $\left \{b_{n};~ n \geq 1 \right \}$ be increasing sequences of 
positive real numbers satisfying (3.1). Then there exist two continuous and increasing functions $\varphi(\cdot)$ and $\psi(\cdot)$
defined on $[0, \infty)$ such that (1.1) holds and 
\begin{equation}
\varphi(0) = \psi(0) = 0, ~\varphi(n) = a_{n}, ~\psi(n) = b_{n}, ~ n \geq 1.
\end{equation}
\end{lemma}

\noindent {\it Proof}~~Let $a_{0} = b_{0} = 0$. Let
\[
\varphi(t) = a_{n-1} + \left(a_{n} - a_{n-1} \right) (t - n + 1), ~n - 1 \leq t < n, ~n \geq 1
\]
and
\[
\psi(t) = b_{n-1} + \left(b_{n} - b_{n-1} \right) (t - n + 1), ~n - 1 \leq t < n, ~n \geq 1.
\]
Clearly, $\varphi(\cdot)$ and $\psi(\cdot)$ are two continuous and increasing functions 
defined on $[0, \infty)$ such that (3.5) holds. We now verify that (1.1) holds with
the chosen $\varphi(\cdot)$ and $\psi(\cdot)$. Note that (3.1) implies that, 
for $n-1 < t < n$ and $n \geq 1$,
\[
\left(\frac{\psi(t)}{\varphi(t)} \right)^{\prime} = 
\frac{\psi^{\prime}(t) \varphi(t) - \psi(t) \varphi^{\prime}(t)}{\varphi^{2}(t)}
= \frac{a_{n-1}b_{n} - b_{n-1}a_{n}}{\varphi^{2}(t)} 
= \frac{a_{n-1}a_{n} \left(\frac{b_{n}}{a_{n}} - \frac{b_{n-1}}{a_{n-1}} \right)}{\varphi^{2}(t)}
\geq 0,  
\]
where $b_{0}/a_{0} \stackrel{\Delta}{=} b_{1}/a_{1}$. Thus (1.1) follows. ~$\Box$

\vskip 0.2cm

\begin{lemma}
{\rm (Corollary 1.3 of Li, Liang, and Rosalsky (2017))}
Let $\{X, X_{n}; ~n \geq 1\}$ be a sequence of i.i.d. 
{\bf B}-valued random variables. Let $\{X_{n}^{\prime};~n \geq 1 \}$ be an
independent copy of $\{X_{n};~n \geq 1 \}$. Write $S_{n} = \sum_{i=1}^{n} X_{i}$,
$S_{n}^{\prime} = \sum_{i=1}^{n} X_{i}^{\prime}$, $n \geq 1$. Let $\{b_{n}; n \geq 1\}$ 
be an increasing sequence of positive real numbers such that 
$\lim_{n \rightarrow \infty} b_{n} = \infty$. Then we have 
\[
\frac{S_{n}- n \mathbb{E}\left(XI\{\|X\| \leq b_{n} \} \right)}{b_{n}}  \rightarrow_{\mathbb{P}} 0
\]
if and only if 
\[
\frac{S_{n} - S_{n}^{\prime}}{b_{n}} \rightarrow_{\mathbb{P}} 0.
\]
\end{lemma}

\vskip 0.3cm

With the preliminaries accounted for, Theorem 3.1 may be proved.

\vskip 0.3cm

\noindent {\it Proof of Theorem 3.1}~~To establish the conclusion of Theorem 3.1,
it suffices to show that, for every sequence $\{X, X_{n}; ~n \geq 1 \}$ 
of i.i.d. {\bf B}-valued random variables, the following three statements are equivalent:
\begin{equation}
\frac{S_{n}- \gamma_{n}}{b_{n}} 
\rightarrow_{\mathbb{P}} 0,
\end{equation}
\begin{equation}
\lim_{n \rightarrow \infty} \mathbb{P}
\left( \frac{\left\|S_{n}- \gamma_{n} \right\|}{b_{n}} > \lambda \right) 
= 0 ~~\mbox{for some constant}~\lambda \in (0, \infty),
\end{equation}
\begin{equation}
\lim_{n \rightarrow \infty}
n \mathbb{P}(\|X\| > b_{n}) = 0.
\end{equation}
Here and below $\gamma_{n} = n \mathbb{E}\left(XI\{\|X\| \leq b_{n} \} \right)$, $n \geq 1$.

Let $\left \{X^{\prime}, X^{\prime}_{n};~ n \geq 1 \right \}$ be an independent copy of
$\left \{X, X_{n};~ n \geq 1 \right \}$ and set $S^{\prime}_{n} = \sum_{i=1}^{n} X^{\prime}_{i}, ~n \geq 1$. 
Since (3.7) obviously follows from (3.6), it suffices to establish the implications ``(3.6) $\Rightarrow$ (3.8)", 
``(3.8) $\Rightarrow$ (3.6)", and ``(3.7) $\Rightarrow$ (3.6)". 

It follows from (3.6) that 
\[
\frac{S_{n} - S_{n}^{\prime}}{b_{n}} \rightarrow_{\mathbb{P}} 0.
\]
and hence by the remarkable L\'{e}vy inequality in a Banach space setting 
(see, e.g., see Proposition 2.3 of Ledoux and Talagrand (1991)), we have that for every $n \geq 1$,
\begin{equation}
\mathbb{P} \left(\frac{\max_{1 \leq i \leq n}\|X_{i} - X_{i}^{\prime} \|}{b_{n}} > t \right) 
\leq 2 
\mathbb{P} \left(\frac{\left \| S_{n} - S_{n}^{\prime} \right \|}{b_{n}} > t \right) \rightarrow 0
~~\mbox{as}~n \rightarrow \infty~~\forall~t \geq 0.
\end{equation}
Note that $\{X - X^{\prime}, X_{n} - X_{n}^{\prime}; ~n \geq 1 \}$ is a sequence of i.i.d. {\bf B}-valued random 
variables. Thus (3.9) implies that
\[
\mathbb{P} \left(\frac{\max_{1 \leq i \leq n}\|X_{i} - X_{i}^{\prime} \|}{b_{n}} > t \right) =
1 - \left(1 - \mathbb{P} \left(\frac{\|X - X^{\prime}\|}{b_{n}} > t \right) \right)^{n}
\rightarrow 0 ~~\mbox{as}~~n \rightarrow \infty ~~\forall~t \geq 0
\]
which is equivalent to
\[
n \mathbb{P} \left(\frac{\|X - X^{\prime}\|}{b_{n}} > t \right) \rightarrow 0 ~~\mbox{as}~~n \rightarrow 
\infty ~~ \forall ~t > 0
\]
and hence (3.8) holds since
\[
\left \{\|X^{\prime} \| \leq b_{n}/2,~\|X\| > b_{n} \right \}
\subseteq \left \{\|X - X^{\prime} \| > b_{n}/2 \right \}
~~\mbox{and}~~\lim_{n \rightarrow \infty} \mathbb{P} \left(\|X^{\prime} \| \leq b_{n}/2 \right) = 1
\]
ensures that for all large $n$,
\[
\mathbb{P} \left (\|X\| > b_{n} \right) \leq 2 \mathbb{P} \left( \|X - X^{\prime} \| > b_{n}/2 \right).
\]
Thus we see that (3.6) implies (3.8).

We now show that (3.8) implies (3.6). Suppose that (3.8) holds. Set
\[
\tilde{X} = \frac{X - X^{\prime}}{2}, ~~\tilde{X}_{n} = \frac{X_{n} - X^{\prime}_{n}}{2}, ~~n \geq 1.
\]
Clearly,
\[
\mathbb{P} (\|\tilde{X} \| > t) \leq \mathbb{P} (\|X\| > t) + \mathbb{P} (\|X^{\prime}\| > t) 
= 2 \mathbb{P} (\|X\| > t) ~~\forall ~t > 0.
\]
Thus $\{\tilde{X}, \tilde{X}_{n}; ~n \geq 1 \}$ is a sequence of i.i.d. 
symmetric {\bf B}-valued random variables such that
\begin{equation}
\lim_{n \rightarrow \infty} n \mathbb{P} \left(\|\tilde{X}\| > b_{n} \right) = 0.
\end{equation}
Since $\{a_{n}; n \geq 1\}$ and $\{b_{n}; n \geq 1\}$ are increasing
sequences of positive real numbers satisfying (3.1), by Lemma 3.1 there exist
two continuous and increasing functions $\varphi(\cdot)$ and $\psi(\cdot)$
defined on $[0, \infty)$ such that both (1.1) and (3.5) hold. Write
\[
Y = \varphi\left(\psi^{-1}(\|\tilde{X}\|)\right) 
\frac{\tilde{X}}{\|\tilde{X}\|}, ~ Y_{n} = \varphi\left(\psi^{-1}(\|\tilde{X}_{n}\|)\right) 
\frac{\tilde{X}_{n}}{\|\tilde{X}_{n}\|}, ~n \geq 1.
\]
It is easy to see that
\[
\mathbb{P}\left(\|Y\| > a_{n} \right) =
\mathbb{P}\left( \left\|\varphi\left(\psi^{-1}(\|\tilde{X}\|)\right) 
\frac{\tilde{X}}{\|\tilde{X}\|} \right \| > \varphi(n) \right) =
\mathbb{P} \left(\|\tilde{X} \| > b_{n} \right), ~n \geq 1.
\]
It thus follows from (3.10) that $\{Y, Y_{n}; ~n \geq 1 \}$ 
is a sequence of i.i.d. symmetric {\bf B}-valued random variables such that
\[
\lim_{n \rightarrow \infty} n \mathbb{P}\left(\|Y\| > a_{n} \right) = 0
\]
and hence by (3.2) 
\begin{equation}
\frac{\sum_{i=1}^{n} Y_{i}}{a_{n}} \rightarrow_{\mathbb{P}} 0.
\end{equation}
By Theorem 1.1 (ii) together with (3.11) and (3.10), we have that for $\epsilon > 0$,
\[
\begin{array}{lll}
\mbox{$\displaystyle
\mathbb{P}\left(\left\|\sum_{i=1}^{n} \tilde{X}_{i} \right\| > \epsilon b_{n} \right)$}
& \leq &
\mbox{$\displaystyle  
4 \mathbb{P} \left(\left\|\sum_{i=1}^{n} \varphi\left(\psi^{-1}(\|\tilde{X}_{i}\|)\right) 
\frac{\tilde{X}_{i}}{\|\tilde{X}_{i}\|} \right\| > \epsilon a_{n} \right) 
+ \sum_{i=1}^{n}\mathbb{P}\left(\|\tilde{X}_{i}\| > b_{n} \right) $}\\
&&\\
&=& 
\mbox{$\displaystyle  
4 \mathbb{P} \left(\left\|\sum_{i=1}^{n} Y_{i} \right\| > \epsilon a_{n} \right) 
+ n \mathbb{P} \left(\|\tilde{X}\| > b_{n} \right)$} \\
&&\\
& \rightarrow &
\mbox{$\displaystyle 0 ~~\mbox{as}~ n \rightarrow \infty $}
\end{array}
\]
and hence
\[
\frac{S_{n} - S_{n}^{\prime}}{2b_{n}} = \frac{\sum_{i=1}^{n} \tilde{X}_{i}}{b_{n}} 
\rightarrow_{\mathbb{P}} 0.
\]
We thus conclude that
\[
\frac{S_{n} - S_{n}^{\prime}}{b_{n}} \rightarrow_{\mathbb{P}} 0.
\]
By Lemma 3.2, (3.6) follows. 

It remains to show that (3.7) implies (3.6). Let $\{X, X_{n}; ~n \geq 1 \}$ 
be a sequence of i.i.d. {\bf B}-valued random variables satisfying (3.7). 
By the L\'{e}vy inequality in a Banach space setting, we have that for all $n \geq 1$,
\[
\mathbb{P} \left( \frac{\max_{1 \leq i \leq n} \left \|X_{i} - X_{i}^{\prime} \right\|}{b_{n}}
> 2 \lambda \right) 
\leq 2 \mathbb{P} \left( \frac{\left \|S_{n} - S_{n}^{\prime} \right\|}{b_{n}}
> 2 \lambda \right) 
\leq 4 \mathbb{P} \left( \frac{\left \|S_{n} \right\|}{b_{n}}
> \lambda \right) ~~\forall~\lambda > 0.
\]
Then it follows from (3.7) that
\[
\lim_{n \rightarrow \infty } 
n \mathbb{P} \left( \frac{\|X - X^{\prime}\|}{b_{n}} > 2 \lambda \right)= 0; ~~\mbox{i.e.,}~
\lim_{n \rightarrow \infty } 
n \mathbb{P} \left( \left\|\frac{X - X^{\prime}}{2 \lambda} \right| > b_{n} \right)= 0 ~~\forall ~\lambda > 0.
\]
That is, (3.8) holds with $X$ replaced by symmetric random variable $(X - X^{\prime})/(2\lambda)$. 
Since (3.6) and (3.8) are equivalent, we conclude that
\[
\frac{\sum_{i=1}^{n} \frac{X_{i} - X_{i}^{\prime}}{2 \lambda}}{b_{n}} \rightarrow_{\mathbb{P}} 0;
~~\mbox{i.e.,}~ \left(\frac{1}{2 \lambda} \right) \frac{S_{n} - S_{n}^{\prime}}{b_{n}} \rightarrow_{\mathbb{P}} 0.
\]
Thus
\[
\frac{S_{n} - S_{n}^{\prime}}{b_{n}} \rightarrow_{\mathbb{P}} 0
\]
which, by Lemma 3.2, implies (3.6). ~$\Box$

\vskip 0.5cm

\noindent
{\bf Acknowledgments}\\

\noindent The research of Deli Li was partially supported by a grant from the Natural Sciences and 
Engineering Research Council of Canada (grant \#: RGPIN-2014-05428)and the research of Han-Ying Liang 
was partially supported by the National Natural Science Foundation of China (grant \#: 11271286).

\vskip 0.5cm

{\bf References}

\begin{enumerate}

\item Chow, Y.S., and Teicher, H. 1997. {\it Probability Theory:
Independence, Interchangeability, Martingales, 3rd ed.}
Springer-Verlag, New York.
    
\item Feller, W. 1946. A limit theoerm for random variables with infinite
moments. Amer. J. Math. {\bf 68}:257-262.
    
\item Hoffmann-J{\o}rgensen, J. 1974. Sums of independent Banach space
valued random variables. {\it Studia Mathematica} {\bf 52}:159-186.

\item Kahane, J.-P. 1985. {\it Some Random Series of Functions, 2nd ed.} Heath Math. Monographs,
1968. Cambridge Univ. Press.

\item Klass, M., and Teicher, H. 1977. Iterated logarithm laws for asymmetric random variables 
barely with or without finite mean. {\it Ann. Probab.} {\bf 5}:861-874.

\item Ledoux, M., and Talagrand, M. 1991. {\it Probability in Banach
Spaces: Isoperimetry and Processes.} Springer-Verlag, Berlin.

\item L\'{e}vy, P. 1937. {\it Th\'{e}orie de L'addition des Variables
Al\'{e}atoires.} Gauthier-Villars, Paris.

\item Li, D., Liang, H.-Y., and Rosalsky, A. 2017. A note on symmetrization procedures 
for the laws of large numbers. {\it Statist. Probab. Lett.}  {\bf 121}: 136-142

\item Marcus, M. B., and Woyczy\'{n}ski, W. A. 1979. Stable measures and
central limit theorems in spaces of stable type. {\it Trans. Amer. Math.
Soc.} {\bf 251}:71-102.

\end{enumerate}

\end{document}